\documentclass[11pt,twoside]{article}

\usepackage{fancyhdr}
\usepackage{epsfig}
\usepackage{latexsym}
\usepackage{microtype}
\usepackage{lastpage} 
\usepackage{hyperref}
\usepackage{graphicx}
\usepackage{lipsum}
\usepackage{amsmath}
\usepackage{amssymb}
\usepackage[dvipsnames*,svgnames]{xcolor}
\usepackage{graphicx}
\usepackage{hyperref}
\usepackage{algorithmic}
\usepackage{tikz}
\usetikzlibrary{arrows}
\usetikzlibrary{positioning}
\usetikzlibrary{decorations}
\usetikzlibrary{decorations.pathmorphing}
\usepackage{multirow}
\usepackage{comment}
\usepackage{caption}
\usepackage{array}
\hypersetup{colorlinks=true,allcolors=DarkBlue,linkcolor=black}

\newcommand{\trshortyear}{20}
\newcommand{\trpapernumber}{02}
\newcommand{\trmonth}{April}
\newcommand{\tryear}{2020}

\DeclareMathOperator*{\argmin}{arg\,min}
\def\D{{\mathcal D}}
\def\E{{\mathcal E}}
\def\H{{\mathcal H}}

\def\S{{\mathcal S}}
\def\T{{\mathcal T}}
\def\R{{\mathbb R}}

\def\Z{{\mathbb Z}}
\def\x{{\mathbf x}}

\newcommand{\TheAuthor}{}
\newcommand{\Author}[1]{\renewcommand{\TheAuthor}{#1}}

\newcommand{\TheTitle}{}
\newcommand{\Title}[1]{\renewcommand{\TheTitle}{#1}}

\Author{E. F. OLARIU, C. FR\u ASINARU}
\Title{Multiple-Depot Vehicle Scheduling Problem Heuristics}
\newcommand\blfootnote[1]{%
	\begingroup
	\renewcommand\thefootnote{}\footnote{#1}%
	\addtocounter{footnote}{-1}%
	\endgroup
}

\begin{document}

	\blfootnote{The publisher does not claim any copyright for the technical reports. The author keeps the full copyright for the paper, and is thus free to transfer the copyright to a publisher if the paper is accepted for publication elsewhere. }   	

\parindent=8mm

\noindent {{\bf \scriptsize Faculty of Computer Science, Alexandru Ioan Cuza University Ia\c si}}

\noindent {{\bf \scriptsize Technical Report TR \trshortyear-\trpapernumber, \trmonth ~ \tryear}}
\vskip -3mm
\noindent\rule{10.2cm}{0.4pt}
\vskip -1mm
\noindent

\vspace{1cm}
\begin{center}
{\Large\bf Multiple-Depot Vehicle Scheduling Problem Heuristics}
\end{center}
\vspace{4mm}

\begin{center}
{\large Emanuel Florentin OLARIU, Cristian FR\u ASINARU}\\
Faculty of Computer Science, Alexandru Ioan Cuza University Ia\c si, \\
General Berthelot 16, 700483 Ias\c i, Romania, \\
Email: {\tt olariu@info.uaic.ro, acf@info.uaic.ro}
\end{center}
\vspace{3ex}

\date{}

\begin{abstract}     
	
The Multiple-Depot Vehicle Scheduling Problem (MDVSP) is very important in the planning process of transport systems. It consists in assigning a set of trips to a set of vehicles in order to minimize a certain total cost. We introduce three fast and reliable heuristics for MDVSP based on a classical integer linear programming formulation and on graph theoretic methods of fixing the infeasible subtours gathered from an integer solution. Extensive experimentations using a large set of benchmark instances show that our heuristics are faster and give good or even better results compared with other existing heuristics.

\smallskip

\noindent
{\bf Keywords:} $ MDVSP$, $ linear \; programming $, $ relaxation \; heuristic $.

\end{abstract}
\section{Introduction}
\label{introduction}

The multiple depot vehicle scheduling problem (MDVSP) is a well-known and important problem combinatorial and optimization problem. MDVSP arises in public transport and trucking industry being part of a large class of problems that includes vehicle routing and scheduling problems. Since these problems are usually NP-hard the exact employed methods cannot solve medium or large MDVSP instances to optimality.

MDVSP aims to assign a set of timetabled tasks (trips) to a set of homogeneous vehicles provided by several depots in order to minimize a given (linear) objective function.


We are given a set of {\it trips} $ T_0, T_1, \ldots, T_{n - 1} $, each trip $ T_i $ having a starting time $ \sigma_i $ and an ending time $ \tau_i $, along with a set of {\it depots}, $ D_0, D_1, \ldots,  D_{m - 1} $, each depot $ D_j $ has a number of $ r_j $ available vehicles. We are given also the time, $ \theta_{ij} $ needed for a vehicle to travel from the end location of trip $ T_i $ to the start location of trip $ T_j $ - these values are useful for deciding if the vehicle performing the trip $ T_i $ can be used after that to perform the trip $ T_j $. An ordered pair of trips $ (T_i, T_j) $ is {\it feasible} if $ \tau_i + \theta_{ij} \le \sigma_j $

A vehicle schedule can be described as an ordered sequence of trips such that any two consecutive trips is a feasible pair. Usually the cost of scheduling includes the sum of the traveling and/or waiting costs between two consecutive trips, the cost of pulling-out the vehicle from the depot to its first trip, and the cost of pulling-in the vehicle from its last trip to the depot.

A solution to MDVSP is an assignment of trips to vehicles that minimizes the sum of costs such that: each trip must be covered, each vehicle schedule must start and end its duty in the same depot, and the number of vehicles available in each depot is not exceeded.

When the number of depots is at least two MDVSP is known to be NP-hard (\cite{bertossi87}). Several approaches have been proposed for this problem; among them: meta-heuristics like neighborhood search and Tabu search (\cite{pepin09}), iterated local search (\cite{laurent09}), and integer linear programming approaches which are the most frequent used methods.

Basically there are three models for MDVSP linear programming formulation: the single-commodity network flow model (introduced in \cite{carpaneto89}), the multi-commodity network flow model and the set partitioning model with side constraints. The multi-commodity model has two different flavors: the classical multi-commodity network flow formulation where the vehicles from different depots are viewed as different commodities (\cite{riberio94}) and the time-space network flow formulation (\cite{kliewer06}, \cite{kulkarni18}). The set partitioning model with side constraints (\cite{riberio94}, \cite{bianco94}, \cite{guedes16}) was derived using Dantzig-Wolfe decomposition.

The paper is organized as follows: section \ref{sect2} describes the linear programming model and two relaxations of this model, section \ref{sect3} describes the ways we fix a solution that contains the so-called {\it infeasible subtours}, section \ref{sect4} contains the numerical results, and section \ref{sect5} concludes the study.

\section{LP model}
\label{sect2}

We use the classical model of single-commodity network flow (\cite{carpaneto89}) because of its reasonable number of variables - in a multi-commodity flow formulation the number of variables is multiplied by the number of depots. The background of this model is the digraph $ G = (V, E) $ (for graph notations and other graph related concepts used in this paper see \cite{diestel00}) where $ V = \{ 0, 1, \ldots, m+ n - 1 \} $. The depots are the vertices in $ V_d = \{ 0, 1, \ldots, m - 1 \} $, while the remaining vertices represent the trips, $ V_t = V \setminus V_d $. 

An arc $ ij $ with $ i, j \in V_t $ exists only if $ (T_{i - m + 1}, T_{j - m + 1}) $ is a pair of feasible trips and has an associated cost $ c_{ij} $ representing the traveling and/or waiting costs between these trips. An arc $ ij $ with $ i \in V_d $ and $ j \in V_d $ (respectively, $ i \in V_t $ and $ j \in V_d $) exists only if the starting (ending) duty of a vehicle from depot $ D_i $ (respectively, $ D_j $) is possible, and its cost $ c_{ij} $ is the incurred cost of starting (ending) the duty from $ D_i $ (respectively to $ D_j $) with the trip $ T_{j - m + 1} $ ($ T_{i - m + 1} $).

The decision variables are $ x_{ij} $ and $ x_{ij} = 1 $ if and only if the arc $ ij $ is used in the optimal solution of the MDVSP. Defining $ r_i $ to be $ 1 $, for $ i \in V_t $ and $ x_{jj} $ to be the the number of unused vehicles in the depot $ D_j $, for $ j \in V_d $, we have the following equivalent linear programming problem:
\begin{align}
&\begin{array}{rcl}
& min & \displaystyle  \left(\sum_{i = 0}^{m +n - 1}\sum_{j = 0}^{m + n - 1} c_{ij}x_{ij} \right)
\end{array}
\label{eqa1}\\
& \qquad \begin{array}{rcl}
\displaystyle \sum_{i = 0}^{m +n - 1} x_{ij} & = & r_j,  0 \le j \le m +n - 1
\end{array}
\label{eqa2}\\
& \qquad \begin{array}{rcl}
\displaystyle \sum_{j = 0}^{m +n - 1} x_{ij} & = & r_i,  0 \le i \le m +n - 1
\end{array}
\label{eqa3}\\
&\qquad \begin{array}{rcl}
\displaystyle \sum_{ij \in E(P)} x_{ij} & \le & |E(P)| - 1, P \in \Pi
\end{array} 
\label{eqa4}\\
&\qquad \begin{array}{rcl}
\displaystyle x_{ij} & \in & \Z_+, 0 \le i, j \le m +n - 1
\end{array} 
\label{eqa5}
\end{align}
where $ \Pi $ is the set of the (inclusion-wise minimal) infeasible paths, that is the paths connecting two different depots. We can linear relax this integer problem by replacing integrality constraints (\ref{eqa5}) with
\begin{align}
&\tag{5'} \begin{array}{rcl}
\displaystyle x_{ij} & \ge & 0, 0 \le i, j \le m +n - 1
\end{array} 
\label{eqn:eqa5prim}
\end{align}

Since $ \Pi $ is very large, even for a reasonable number of depots, the constraints \eqref{eqa4} could be difficult be difficult to implement, since using an enumerative technique for finding $ \Pi $ would be very costly. Our approaches are based on relaxing by reducing the number of constraints of type \eqref{eqa4} and, than, apply a technique of repairing the integer solutions.

\subsection{First relaxation - Circulation}
\label{subsect2_1}

The first relaxation of our problem is (\ref{eqa1} - \ref{eqa3}), \eqref{eqa5} which can be solved efficiently by modeling it as a \textit{minimum-cost circulation problem} \cite{ahuja93}, a generalization of the minimum-cost flow problem with node capacities and lower bounds on the edges.

For a given transportation network $R=(G, c, l, u, a)$, where $G=(V,E)$ is a digraph, $c:V \rightarrow \R$ represents the \textit{capacities} of the nodes, $l,u : E \rightarrow \R_+$ represent \textit{lower} and \textit{upper} bounds for the edges and $a: E \rightarrow \R$ is a \textit{cost} function interpreted as the cost of "sending" an unit of flow on a specific arc. 
If $ c(v) > 0 $, $ v $  is called a \textit{supply} node, if $ c(v)  < 0 $, $ v $ is a \textit{demand} node, the remaining nodes being \textit{transit} nodes.
The sum of all node capacities must be zero, meaning that the supply must equal the demand. 
Let $S=\{v \in V \;|\; b(v) > 0 \}$ be the \textit{source} nodes and $T=\{v \in V \;|\; b(v) < 0 \}$ be the \textit{sink} nodes.

The minimum-cost circulation problem is to find a feasible circulation (a function that satisfies equations (\ref{eqa7} - \ref{eqa8})) that minimizes the total cost:
\begin{align}
&\begin{array}{rcl}
& min & \displaystyle  \left( \sum_{ij \in E} c_{ij} x_{ij} \right)
\end{array}
\label{eqa6}\\
&\begin{array}{rcl}
\displaystyle l_{ij} \le x_{ij} \le u_{ij}, \forall ij \in E
\end{array}
\label{eqa7}\\
&\begin{array}{rcl}
\displaystyle \sum_{ij \in E}  x_{ij} - \sum_{ji \in E} x_{ji} = b_i , \forall i \in V
\end{array}
\label{eqa8}
\end{align}

The problem can be solved in a pseudo-polynomial time $O(nU(m+n)\log{n})$ using the successive shortest path algorithm with capacity scaling \cite{ahuja93}, where $n=|V|$, $m=|E|$, and $U$ is the maximum edge capacity.

In order to represent our scheduling problem as a minimum-cost circulation problem we define the following transportation network $R=(G', b, l, u, c')$, based on the initial digraph $ G $ and cost matrix $c$.

\begin{itemize}
\item For each depot $ D_i $, add two nodes in $V(G')$, representing a source and a sink. The capacities of these nodes are $ c(i) = r_i $ (the supply) and $ c(i') = - r_i $ (the demand).

\item For each trip $ T_{j - m+ 1} $, add two nodes $ j^- $ and $ j^+ $ in $ G' $ with zero capacity (they will be transit nodes), and the arc $ j^- j^+ $ in $ G' $, having the cost $0$ (we are transforming the trip nodes in $ G $ into arcs in $ G' $); both the lower and the upper bounds of these arcs are set to $1$, as we are looking for a solution that saturates all trips.

\item For each arc $ ij $ connecting the depot $ D_i $ to a trip $ T_{j - m + 1} $ in $ G $, add the arc $ ij^- $ in $ G' $, having the cost $ c_{ij} $, the lower bound $ 0 $, and the upper bound $ 1 $.

\item For each arc $ ji $ connecting the trip $ T_{j - m + 1}$ to a depot $ D_i $ in $ G $, add the arc $j^+i$ in $ G' $, having the cost $c_{ji}$, the lower bound $ 0 $, and the upper bound $ 1 $.

\item For each arc $ ij $ connecting two trips $ T_{i - m + 1} $ and $ T_{j - m + 1} $ in $ G $, add the edge $i^-j^+$ in $ G' $, having the cost $c_{ij}$, the lower bound $ 0 $, and the upper bound $ 1 $.

\item For each depot $ D_i $ in $ G $, add the arc $ ii' $ in $ G' $, having the cost $ 0 $, the lower bound $ 0 $, and the upper bound $ r_i $ (this arc is necessary when some vehicles in the depot $ D_i $ are not used). 

\end{itemize}

It is straightforward to prove that a solution for the minimum-cost circulation problem defined for the network $R=(G', b, l, u, c') $ represents an optimal solution for the relaxation (\ref{eqa1} - \ref{eqa3}), \eqref{eqa5}.
The total cost is the same for both problems, the supply/demand constraints ensure that the number of vehicles used from a depot is not exceeded and the lower bound constraints imposed for the arcs $j^-j^+$ ensure that all trips are saturated.

This model is a guarantee that the relaxed MDVSP problem can be solved in an efficient manner. From a practical point of view, the polynomial complexity guarantee usually translates in an easy resolution when it comes to dedicated MIP solvers, such as Gurobi.

\subsection{Second relaxation}
\label{subsect2_2}

Our second method of relaxing the original problem consists in replacing $ \Pi $ by a smaller set of infeasible subtours (paths that link different depots) by adding one by one constraints of type \eqref{eqa4} and re-optimizing until the new problem has the same optimum as (\ref{eqa1} - \ref{eqa4}), \eqref{eqn:eqa5prim}. This procedure of re-optimizing is based on the interpretation of graph theoretic properties of a fractional solution.

At a certain step during the algorithm we have a particular set of infeasible subtours $ \Pi' $ and \eqref{eqa4} is replaced in the current problem by
\begin{align}
\tag{4'}
&\begin{array}{rcl}
\displaystyle \sum_{ij \in E(P)} x_{ij} & \le & |E(P)| - 1, P \in \Pi'
\end{array} 
\label{eqn:eqa4prim}
\end{align}
Consider $ \x^* = (x^*_{ij}) $, a solution to problem (\ref{eqa1} - \ref{eqa3}), \eqref{eqn:eqa4prim}, \eqref{eqn:eqa5prim}, and define a weight on the edges of the underlying digraph: $ \alpha_{ij} = 1 - x^*_{ij} $, for all arcs $ ij $. \eqref{eqa4} is equivalent with
\begin{align}
\tag{4''}
&\begin{array}{rcl}
\alpha(P) \ge 1, P \in \Pi
\end{array}
\label{eqn:eqa4sec}
\end{align}
since
\begin{align}
\notag
& \sum_{ij \in E(P)} x^*_{ij} \le |E(P)| - 1 \Leftrightarrow \sum_{ij \in E(P)} (1 - x^*_{ij}) \ge 1 \Leftrightarrow \sum_{ij \in E(P)} \alpha_{ij} \ge 1.
 \end{align}
Hence, $ \x^*$ is an optimum solution to (\ref{eqa1} - \ref{eqa4}), \eqref{eqn:eqa5prim} if and only if the underlying digraph doesn't contain paths between different depots of sub-unitary weight. We will test this by using an algorithm for finding shortest paths in a weighted graph, like Floyd-Warshall or Bellman-Ford-Moore.

Therefore the first step is to relax the problem (\ref{eqa1} - \ref{eqa4}), \eqref{eqn:eqa5prim} to (\ref{eqa1} - \ref{eqa3}), (\ref{eqn:eqa4prim} - \ref{eqn:eqa5prim}) for a certain known set of infeasible paths $ \D' $ such that the two problems have the same optimum. The process of building this problem is given below.
 \begin{algorithmic}
\STATE $ \Pi' \gets \varnothing $;
\STATE solve problem (\ref{eqa1} - \ref{eqa3}), (\ref{eqn:eqa4prim} - \ref{eqn:eqa5prim}) and let $ \x^* $ be an optimum solution;
\WHILE{(there exists a path $ D $ with $ \alpha(D) < 1 $)}
	\STATE add $ D $ to $\D' $;
	\STATE solve the problem (\ref{eqa1} - \ref{eqa3}), (\ref{eqn:eqa4prim} - \ref{eqn:eqa5prim}) and let $ \x^* $ be an optimum solution;
\ENDWHILE
\STATE return $ \x^* $.

\end{algorithmic}

The aim of the above procedure is to build a problem that has a larger (but known) set of feasible solutions but the same optimum with (\ref{eqa1} - \ref{eqa4}), \eqref{eqn:eqa5prim}
 
\subsection{Column generation perspective for the second relaxation}
\label{subsect2_3}

Writing the original relaxed problem as a maximum one means to replace \eqref{eqa1} by  \eqref{eqa1prim} (ignoring the minus in front of {\it max})
\begin{align}
\tag{1'}
&\begin{array}{rcl}
& max & \displaystyle  \left( \sum_{i = 0}^{m +n - 1}\sum_{j = 0}^{m + n - 1} -c_{ij}x_{ij} \right) 
\end{array}
\label{eqa1prim}
\end{align}
The dual of the problem \eqref{eqa1prim}, \eqref{eqa2} - \eqref{eqa4}, \eqref{eqn:eqa5prim} is

\begin{align}
&\begin{array}{rcl}
& min & \displaystyle \left( \sum_{j = 0}^{m +n - 1}r_jy_j + \sum_{i = 0}^{m + n - 1} r_iz_i + \sum_{P \in \Pi} (|E(P)| - 1) u_P \right) 
\end{array}
\label{eqad1}\\
& \begin{array}{rcl}
\displaystyle z_i + y_j + \sum_{P \in \Pi: ij \in E(P)} u_P \ge -c_{ij}, \forall ij \in E(G)
\end{array}
\label{eqad2}\\
&\begin{array}{rcl}
\displaystyle u_P \ge 0, \forall P \in \Pi
\end{array}
\label{eqad3}
\end{align}
We can replace \eqref{eqad2} by \eqref{eqad2prim}, and  \eqref{eqad3} by \eqref{eqad3prim}

\begin{align}
\tag{10'}
&\begin{array}{rcl}
\displaystyle z_i + y_j + \sum_{P \in \Pi: ij \in E(P)} u_P -v_{ij} = -c_{ij}, \forall ij \in E(G)
\end{array}
\label{eqad2prim}\\
& \tag{11'}
\begin{array}{rcl}
\displaystyle v_{ij} \ge 0, \forall ij \in E(G), u_P \ge 0, \forall P \in \Pi
\end{array}
\label{eqad3prim}
\end{align}
and get the dual in equations form: \eqref{eqad1} \eqref{eqad2prim},  \eqref{eqad3prim}. An initial feasible basic solution to this problem could be $ v_{ij} = c_{ij} $, $ \forall ij \in E(G) $. Now, this dual problem has a very large number of variables and we can solve it by using the column generation method (see for example \cite{desrosiers05}, \cite{desrosiers04}). We start with a small set of variables (that contains a feasible basis)  - this is the restricted master problem - and in each step  - by solving the corresponding sub-problem - find a variable with the minimum negative reduced cost that would be added to the current problem. When such variables doesn't exists we have an optimum solution to the dual problem.

In our case the sub-problem would be
\begin{align}
\tag{10'}
&\begin{array}{rcl}
\displaystyle \argmin_{P \in \Pi}\left( |E(P)| - 1 - \sum_{ij \in E(P)} x_{ij} \right) < 0
\end{array}
\label{eqadsubpr}
\end{align}
This holds because the variables $ z_i $ and $ y_j $ cannot enter the basis - they cannot have negative reduced cost based on \eqref{eqa2} and \eqref{eqa3}.

Adding a new variable (column) to the dual problem means adding a new constraint to the primal, i. e., finding a new path between different depots of $ \alpha $ sub-unitary cost. Hence solving the problem with the algorithm from the subsection \ref{subsect2_2} is equivalent with solving the dual using the column generation.

\section{Building feasible solutions - Path repairing methods}
\label{sect3}

The next step is to solve one of the two problems \eqref{eqa1} - \eqref{eqa3}, \eqref{eqa5} or \eqref{eqa1} - \eqref{eqa3}, \eqref{eqn:eqa4prim}, \eqref{eqa5} which are integer linear programming problems. Because of their reasonable size these problems can be easily solved with existing LP solvers - for the first one there exist even combinatorial algorithms (see subsection \ref{subsect2_1}). Solutions to this problems may not be feasible for  \eqref{eqa1} - \eqref{eqa5}, thus a process of fixing the infeasible subtours (that is a path that starts in a depot and ends in a different depot) must follow. 

In the remaining of this section we suppose that we have a solution, $ \x^* $, to one of the above two integer linear programming problems.

First we define an auxiliary digraph $ H = (V_d, A_d) $, where $ ij \in A_d $ if and only if there is an infeasible subtour between the depots $ i $ and $ j $; we add also a weight on this arc $ w_{ij} $ which represents the number of such subtours. By inspecting $ \x^* $ we extract the infeasible tours and memorize them, build $ H $ and the weight $ w $. 

\subsection{Repairing one subtour}

Suppose that we have the infeasible subtour (described by its arcs):
\begin{align}
\notag
P:  it_1, t_1t_2, \ldots, t_{p - 1}t_p, t_pj, \; i, j \in V_d, 
\end{align}
where $ i, j \in V_d $ and $ t_h \in V_t, \forall h = \overline{1, p} $; this subtour can be replaced by
\begin{align}
\notag
P':  it_1, t_1t_2, \ldots, t_{p - 1}t_p, t_pi
\end{align}
with cost penalty $ c'(P) = \displaystyle c\left(t_pi \right) - c\left(t_pj \right) $, or by 
\begin{align}
\notag
P": jt_1, t_1t_2, \ldots, t_{p - 1}t_p, t_pj,
\end{align}
with cost penalty $ c''(P) = \displaystyle c\left(jt_1 \right) - c\left(it_1 \right) $.

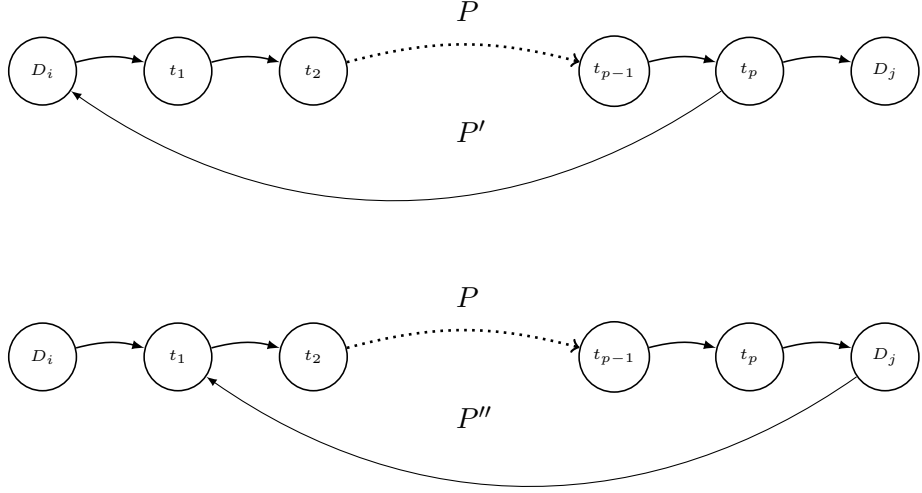
\begin{figure}

\begin {center}

\begin {tikzpicture}[-latex ,auto ,node distance =1.8cm ,on grid , semithick , state/.style ={ circle ,top color =white , draw,black , text=black , minimum width = 0.9 cm}]
\node[state]  (Di)    {\tiny$ D_i $};
\node[state] (t11) [right= of Di] {\tiny$ t_1 $};
\node[state] (t12) [right=of t11] {\tiny$ t_2 $};
\node[state] (t1pp) [right= 4cm of t12] {\tiny $ t_{p - 1} $};
\node[state] (t1p) [right= of t1pp] {\tiny$ t_p $};
\node[state] [right =of t1p] (Dj)  {\tiny$ D_j $};
\path (Di) edge [bend right = -15] node[below =0.15 cm] {} (t11);
\path (t11) edge [bend right = -15] node[below =0.15 cm] {} (t12);
\draw[->,dotted,line width=1pt] (t12) to  [bend right=-15] (t1pp);	
\path (t1pp) edge [bend right = -15] node[below =0.15 cm] {} (t1p);
\path (t1p) edge [bend right = -15] node[below =0.15 cm] {} (Dj);
\path (t1p) edge [thin, bend right = -35] node[below =0.15 cm] {} (Di);
\node[text width=3cm] at (7,0.8) {$ P $};
\node[text width=3cm] at (7,-0.8) {$ P' $};
\end{tikzpicture}

\vspace{0.5cm}

\begin {tikzpicture}[-latex ,auto ,node distance =1.8cm ,on grid , semithick , state/.style ={ circle ,top color =white , draw,black , text=black , minimum width = 0.9 cm}]
\node[state]  (Di)    {\tiny$ D_i $};
\node[state] (t11) [right= of Di] {\tiny$ t_1 $};
\node[state] (t12) [right=of t11] {\tiny$ t_2 $};
\node[state] (t1pp) [right= 4cm of t12] {\tiny $ t_{p - 1} $};
\node[state] (t1p) [right= of t1pp] {\tiny$ t_p $};
\node[state] [right =of t1p] (Dj)  {\tiny$ D_j $};
\path (Di) edge [bend right = -15] node[below =0.15 cm] {} (t11);
\path (t11) edge [bend right = -15] node[below =0.15 cm] {} (t12);
\draw[->,dotted,line width=1pt] (t12) to  [bend right=-15] (t1pp);	
\path (t1pp) edge [bend right = -15] node[below =0.15 cm] {} (t1p);
\path (t1p) edge [bend right = -15] node[below =0.15 cm] {} (Dj);
\path (Dj) edge [thin, bend right = -35] node[below =0.15 cm] {} (t11);
\node[text width=3cm] at (7,0.8) {$ P $};
\node[text width=3cm] at (7,-0.8) {$ P'' $};
\end{tikzpicture}

\end{center}
\caption{Repairing one subtour.}
\label{fig:repair1}
\end{figure}

\subsection{Repairing a pair of subtours}

Suppose now that we have two infeasible subtours
\begin{align}
\notag
P_1: it^1_1, t^1_1t^1_2, \ldots, t^1_{p - 1}t^1_p, t^1_pj \mbox{ and } P_2:  jt^2_1, t^2_1t^2_2, \ldots, t^2_{q - 1}t^2_q, t^2_qi,
\end{align}
where $ i, j \in V_d $ and  $ t^1_h, t^2_k \in V_t, \forall h = \overline{1, p}, k = \overline{1, q} $. Such a pair will be called {\it compatible}.

If we can find a pair $ (h, k) $, $ 1 \le h \le p $ and $ 1 \le k \le q $ such that $ t^1_ht^2_k, t^2_{k - 1}t^1_{h + 1} \in E(G) $, then we can replace the above pair of infeasible but compatible subtours by the following pair of feasible subtours
\begin{align}
\notag
& P_1': it^1_1, t^1_1t^1_2, \ldots, t^1_{h - 1}t^1_h, t^1_ht^2_k, t^2_kt^2_{k + 1}, \ldots, t^2_{q - 1}t^2_q, t^2_qi\\
\notag
& P_2':  jt^2_1, t^2_1t^2_2, \ldots, t^2_{k - 2}t^2_{k - 1}, t^2_{k - 1} t^1_{h + 1}, t^1_{h + 1}t^1_{h + 2}, \ldots,  t^1_{p - 1}t^1_p, t^1_pj. 
\end{align}
The cost penalty is
\begin{align}
\notag
\displaystyle \gamma_{hk} = c\left(t^1_ht^2_k \right) + c\left(t^2_{k - 1}t^1_{h - 1} \right) - c\left(t^1_it^1_{i + 1} \right) - c\left(t^2_{k - 1}t^2_k \right).
\end{align}
We will choose the pair $ (h, k) $ for which the cost penalty is minimum, i. e.,
\begin{align}
\notag
& c(P_1, P_2) = \min_{1 \le h \le p, 1 \le k \le q}{\gamma_{hk}}.
\end{align} 

\begin{figure}


\begin {center}
\begin {tikzpicture}[-latex ,auto ,node distance =1.8cm ,on grid , semithick , state/.style ={ circle ,top color =white , draw,black , text=black , minimum width = 1 cm}]
\node[state]  (Di)  {\tiny$ D_i $};
\node[state] (t11)  [above right = 2.5 cm of Di] {\tiny$ t^1_1 $};

\node[state] (t1h) [right= 3 cm of t11] {\tiny$ t^1_h $};
\node[state] (t1hh) [right= of t1h] {\tiny $ t^1_{h+ 1} $};
\node[state] (t1p) [right= 3cm of t1hh] {\tiny$ t^1_p $};
\node[state] [right =of t1p] (Dj)  {\tiny$ D_j $};
\node[state] (t21)  [below left = 2.5 cm of Dj] {\tiny$ t^2_1 $};
\node[state] (t2kk) [left= 3 cm of t21] {\tiny$ t^2_{k - 1} $};
\node[state] (t2k) [left= of t2kk] {\tiny $ t^2_{k} $};
\node[state] (t2q) [left= 3 cm of t2k] {\tiny$ t^2_q $};
\path (Di) edge [bend right = -15] node[below =0.15 cm] {} (t11);
\draw[->,dotted,line width=1pt] (t11) to  [bend right=-15] (t1h);	
\path (t1h) edge [bend right = -15] node[below =0.15 cm] {} (t1hh);
\draw[->,dotted,line width=1pt] (t1hh) to  [bend right=-15] (t1p);	
\path (t1p) edge [bend right = -15] node[below =0.15 cm] {} (Dj);
\node[text width=3cm] at (7,-0.7) {$ P_2 $};

\path (Dj) edge [bend right = -15] node[below =0.15 cm] {} (t21);
\draw[->,dotted,line width=1pt] (t21) to  [bend right=-15] (t2kk);	
\path (t2kk) edge [bend right = -15] node[below =0.15 cm] {} (t2k);
\draw[->, dotted,line width=1pt] (t2k) to  [bend right=-15] (t2q);
\path (t2q) edge [bend right = -15] node[below =0.15 cm] {} (Di);
\node[text width=3cm] at (7, 2.4) {$ P_1 $};

\path (t1h) edge [thin, bend right = -15] node[below =0.15 cm] {} (t2k);
\path (t2kk) edge [thin, bend right = -15] node[below =0.15 cm] {} (t1hh);

\node[text width=3cm] at (4.5,0.8) {$ P'_1 $};
\node[text width=3cm] at (9.5, 0.8) {$ P'_2 $};

\end{tikzpicture}

\end{center}
\caption{Repairing two subtours.}
\label{fig:repair2}
\end{figure}
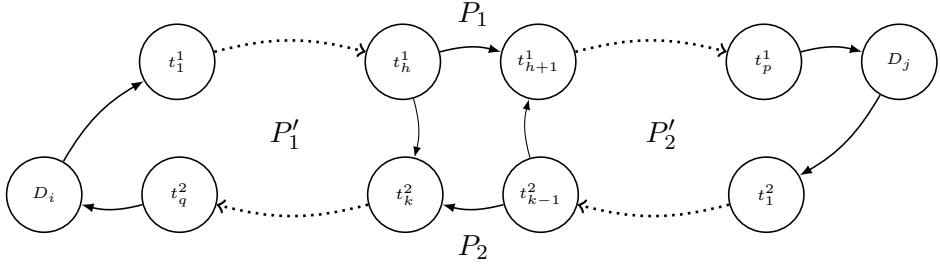

\subsection{Fixing all subtours}

The first way of repairing all subtours is to iterate the first method for all infeasible tours. A better way is to use both methods above: we match infeasible subtours or choose to repair a subtours by itself.

Since for any specific infeasible subtour $ P_1 $ there may be more than one compatible subtour $ P_2 $, we must match infeasible subtours in a manner that minimizes the overall penalty of the repairing. To implement this method we create a bipartite graph $ \H = (\S, \T; \E) $ that models the compatibility relation between infeasible subtours.
The partition classes are
\begin{align}
\notag
& \S = \{ P \: : \: P \mbox{ is an infeasible subtour with respect to } \x^* \} \\
\notag
& \T = \{ P' \: : \: P \mbox{ is an infeasible subtour with respect to } \x^* \}.
\end{align}
The set of edges is
\begin{align}
\notag
\E  = \{ P_1P_2' \: : \: (P_1, P_2) \mbox{ is a compatible pair of infeasible subtours} \} \cup\\  \cup \{ PP' \: : \: P \mbox{ is an infeasible subtour} \}
\end{align}
We define on these edges two weight functions $ \gamma', \gamma'' : \E \to \R $: $ \gamma'(P_1P_2') = \gamma''(P_1P_2') = c(P_1, P_2) $, if $ P_1 $ and $ P_2 $ are compatible infeasible subtours, and $ \gamma'(PP') = c'(P) $, $ \gamma''(PP') = c''(P) $ if $ P $ is an infeasible subtour.

Since $ \H $ has perfect matchings (due to the edges of the form $ PP' $), we can find two $ \gamma' $ - and $ \gamma'' $ - minimum weight perfect matchings using the Kuhn-Munkres (Hungarian) algorithm in $ O(|\S|^3)$ time complexity. We choose the matching having the smaller weight for repairing the solution $ \x^* $; applying this method will fix all the infeasible subtours. In our numerical experiments the number of infeasible subtours is quite small, hence this method will work fast in practice.

In this way the method of repairing just one subtour is used by this second method that fixes all the infeasible tours. Iterating the first method cannot improve the result of the second since the latter already takes account the edges $ PP' $ (for all infeasible subtours $ P $).

\subsection{Heuristics}

We developed three heuristics based on these linear integer problems and fixing subtours method. The first heuristic($ H_1 $), after solving (\ref{eqa6} - \ref{eqa8}) (or, equivalently, (\ref{eqa1} - \ref{eqa3}), \eqref{eqa5}  ), fixes all the infeasible subtours by finding the minimum weight matching in the corresponding bipartite graph. The second heuristic ($ H_2 $) requires a given number (a pool) of integer solutions and then builds the bipartite graphs for all of these solutions, fixes the subtours and chooses the best repaired solution. 

The third heuristic ($ H_3 $) first builds the set of infeasible subtours $ \Pi' $ by repeatedly finding paths of sub-unitary weight in the subjacent digraph, then solves the ILP (\ref{eqa1} - \ref{eqa3}), \eqref{eqn:eqa4prim}, \eqref{eqa5}; the resulting solution is then fixed using the above method. 

\section{Computational results}
\label{sect4}
In this section we describe the results of our numerical experiments using the heuristics described above. Our test bed is composed of the instances used in (\cite{pepin09}) (see the Huisman's website \url{https://personal.eur.nl/huisman/instances.htm}) and most of the instances in (\cite{kulkarni18}) (see (\cite{kulkarni19})), all generated with the method from (\cite{carpaneto89}). 

All computational results from below were obtained using an Intel(R) Core (TM) i5-7500 CPU \@ 3.40GHz computer with 8GB RAM, under Ubuntu 18.04.4 LTS.

The linear programming problems were solved using Gurobi 9.0 under an Academic License.

The experiments were performed using seventy different  instances. Results are very close to the best known results; for instances with at most $ 1500 $ trips and at most $ 8 $ depots all the results but three are within $ 1\% $ relative error from the best known corresponding solutions. The fastest of the three heuristics is, as expected, $ H_1 $; the best results are obtained with $ H_2 $ or $ H_3 $. The differences among the three heuristics are very small.

The results reported in the tables below use the following performance measure:
\begin{align}
\notag
\mbox{percent error (\%)} = 100 \cdot \frac{z - z^0}{z^0},
\end{align}
where $ z^0 $ the the best known objective function value and $ z $ is the heuristic's obtained value.

Tables \ref{table1} and \ref{table2} show the best known objective function values, the objective values for our heuristics, the CPU times spent, and the percent errors. For the larger instances (table \ref{table2}) the second heuristic ($ H_2 $) proved to be too time costly, $ H_1 $ remains the faster and $ H_3 $ gives the best results in terms of gap.

\begin{table}[h]
\caption{Solutions and computational times for different heuristics.}
\label{table1}
\footnotesize
\begin{tabular}{>{\centering}m{0.10\textwidth}>{\centering}m{0.10\textwidth}>{\centering}m{0.10\textwidth}>{\centering}m{0.10\textwidth}>{\centering}m{0.10\textwidth}>{\centering}m{0.03\textwidth}>{\centering}m{0.03\textwidth}>{\centering}m{0.03\textwidth}>{\centering}m{0.04\textwidth}>{\centering}m{0.04\textwidth}c}
\hline
   \multirow{2}{*}{Instance} & Best & \multicolumn{3}{c}{Heuristics solutions} & \multicolumn{3}{c}{CPU time (s)} &  \multicolumn{3}{c}{Percent error (\%)} \\
\cline{3-11}
& solution & $ H_1 $ & $ H_2 $ & $ H_3 $ &  $ H_1 $ & $ H_2 $ & $ H_3 $ &  $ H_1 $ & $ H_2 $ & $ H_3 $ \\
\hline
m4n500s0 & 1,289,114 & 1,296,409 & {\bf 1,295,671}  & 1,295,678 &  0.7 & 7.9 & 6.5 &  0.56\% & 0.50\% & 0.50\% \\
m4n500s1 & 1,241,618 & 1,247,438 & {\bf 1,246,655}  & 1,247,173 &  0.6 & 11.8 & 5.4 &  0.46\% & 0.40\% & 0.44\% \\
m4n500s2 & 1,283,811 & 1,292,079 & 1,291,745  & {\bf 1,290,891} &  0.5 & 17.4 & 5.6 &  0.64\% & 0.61\% & 0.55\% \\
m4n500s3 & 1,258,634 & 1,263,624 & {\bf 1,263,045}  & 1,264,473 &  0.4 & 10.2 & 4.7 &  0.39\% & 0.35\% & 0.46\% \\
m4n500s4 & 1,317,077 & 1,322,535 & 1,322,306  & {\bf 1,321,138} &  0.5 & 11.7 & 4.9 &  0.41\% & 0.39\% & 0.30\% \\
\hline
m4n1000s0 & 2,516,247 & 2,528,728 &{\bf  2,527,966}  & 2,528,299  & 2.8 & 41.4 & 48.7 &  0.49\% & 0.46\% & 0.47\% \\
m4n1000s1 & 2,413,393 & 2,421,735 & 2,421,735  & {\bf 2,420,440} &  2.5 & 37.4 & 29.2 &  0.34\% & 0.34\% & 0.29\% \\
m4n1000s2 & 2,452,905 & 2,461,985 & 2,461,787  & {\bf 2,461,347} &  2.7 & 38.4 & 17.9 &  0.37\% & 0.36\% & 0.34\% \\
m4n1000s3 & 2,490,812 & 2,498,319 & {\bf 2,498,046}  & 2,498,423 &  2.5 & 65.3 & 28.8 &  0.30\% & 0.29\% & 0.30\% \\
m4n1000s4 & 2,519,191 & 2,525,357 & 2,525,004  & {\bf 2,524,898} &  2.6 & 45.6 & 13.5 &  0.24\% & 0.23\% & 0.22\% \\
\hline
m4n1500s0 & 3,830,912 & 3,847,046 & 3,846,785 & {\bf 3,846,761} & 6.5 & 151.6 & 136.6 &  0.41\% & 0.41\% & 0.41\% \\
m4n1500s1 & 3,559,176 & 3,566,055 & 3,565,962 & {\bf 3,564,918} &  7.0 & 177.1 & 56.5 &  0.19\% & 0.19\% & 0.16\% \\
m4n1500s2 & 3,649,757 & 3,662,319 & 3,661,323 & {\bf 3,661,344} &  7.2 & 232.1 & 115.4 & 0.34\% & 0.31\% & 0.31\% \\
m4n1500s3 & 3,406,815 & 3,419,905 & 3,419,810 & {\bf 3,417,225} &  5.9 & 289.2 & 118.7 &  0.38\% & 0.38\% & 0.30\% \\
m4n1500s4 & 3,567,122 & 3,583,176 & 3,582,852 & {\bf 3,581,059} &  6.3 & 208.8 & 122.5 &  0.45\% & 0.44\% & 0.39\% \\
\hline
m8n500s0 & 1,292,411 & 1,304,837 & 1,302,517  & {\bf 1,301,395} &  0.5 & 9.4 & 8.1 &  0.96\% & 0.78\% & 0.69\% \\
m8n500s1 & 1,276,919 & 1,289,875 & 1,288,006  & {\bf 1,289,407} &  0.6 & 42.4 & 5.4 &  1.01\% & 0.86\% & 0.97\% \\
m8n500s2 & 1,304,251 & 1,316,965 & 1,316,108  & {\bf 1,313,993} &  0.5 & 18.5 & 5.8 &  0.97\% & 0.90\% & 0.74\% \\
m8n500s3 & 1,277,838 & {\bf 1,290,397} & {\bf 1,290,397} & 1,290,852 &  0.5 & 9.4 & 6.4 &  0.98\% & 0.98\% & 1.01\% \\
m8n500s4 & 1,276,010 & 1,289,435 & {\bf 1,287,919}  & 1,288,606 &  0.6 & 11.8 & 7.1 &  1.05\% & 0.93\% & 0.98\% \\
\hline
m8n1000s0 & 2,422,112 & 2,441,490 & {\bf 2,439,817}  & 2,439,893 &  2.7 & 143.5 & 57.6 &  0.80\% & 0.73\% & 0.73\% \\
m8n1000s1 & 2,524,293 & 2,542,668 & {\bf 2,542,668}  & 2,545,417 &  3.0 & 31.3 & 33.2 &  0.72\% & 0.72\% & 0.83\% \\
m8n1000s2 & 2,556,313 & 2,581,639 & 2,580,507  & {\bf 2,579,511} &  2.6 & 182.0 & 45.4 &  0.99\% &  0.94\% & 0.90\% \\
m8n1000s3 & 2,478,393 & 2,499,109 & 2,495,968  & {\bf 2,494,389} &  2.8 & 217.5 & 38.0 &  0.83\% & 0.70\% & 0.64\% \\
m8n1000s4 & 2,498,388 & 2,518,121 & 2,517,631  & {\bf 2,516,357} &  2.9 & 34.3 & 42.1 &  0.79\% & 0.77\% & 0.71\% \\
\hline
m8n1500s0 & 3,500,160 & 3,527,083 & {\bf 3,527,083}  & 3,530,381 &  5.8 & 114.6 & 106.9 &  0.76\% & 0.76\% & 0.86\% \\
m8n1500s1 & 3,802,650 & 3,821,483 & 3,819,634  & {\bf 3,818,617} &  6.0 & 335.2 & 393.9 &  0.49\% & 0.44\% & 0.42\% \\
m8n1500s2 & 3,605,094 & 3,640,171 & {\bf 3,635,622} & 3,636,799 & 5.5 & 268.3 & 155.4 &  0.97\% & 0.84\% & 0.87\% \\
m8n1500s3 & 3,515,802 & 3,537,090 & {\bf 3,536,906} & 3,536,931 & 5.8 & 283.8 & 139.1 &  0.60\% & 0.60\% & 0.60\% \\
m8n1500s4 & 3,704,953 & 3,733,572 & 3,733,572 & {\bf 3,730,221} & 5.5 & 101.5 & 120.7 &  0.77\% & 0.77\% & 0.68\% \\
    \end{tabular}
\end{table}

\begin{table}[h]
\caption{Solutions and computational times for different heuristics.}
\label{table2}
\footnotesize
\begin{tabular}{>{\centering}m{0.11\textwidth}>{\centering}m{0.10\textwidth}>{\centering}m{0.10\textwidth}>{\centering}m{0.08\textwidth}>{\centering}m{0.08\textwidth}>{\centering}m{0.08\textwidth}>{\centering}m{0.08\textwidth}c}
\hline
   \multirow{2}{*}{Instance} & Best & \multicolumn{2}{c}{Heuristics solutions} & \multicolumn{2}{c}{CPU time (s)} &  \multicolumn{2}{c}{Percent error (\%)} \\
\cline{3-8}
& solution & $ H_1 $  & $ H_3 $ &  $ H_1 $ & $ H_3 $ &  $ H_1 $ & $ H_3 $ \\
\hline
m8n2000s0 & 4,916,810 & 4,975,718 & {\bf 4,962,626} & 9.8 & 402.9 & 1.19\% & 0.93\% \\
m8n2000s1 & 4,769,442 & 4,819,440 & {\bf 4,813,103} & 9.4 & 807.6 & 1.04\% & 0.91\% \\
m8n2000s2 & 4,897,886 & 4,948,430 & {\bf 4,938,756} & 8.9 & 474.7 & 1.03\% & 0.83\% \\
m8n2000s3 & 5,171,924 & 5,231,090 & {\bf 5,220,119} & 9.0 & 513.5 & 1.14\% & 0.93\% \\
m8n2000s4 & 4,761,862 & 4,808,420 &  {\bf 4,802,721} & 9.1 & 469.5 & 0.97\% & 0.85\% \\
\hline
m8n2500s0 & 5,911,824 & 5,981,468 & {\bf 5,961,055} &  14.7 & 879.1 & 1.17\% & 0.83\% \\
m8n2500s1 & 6,296,870 & 6,363,706 & {\bf 6,357,577} &  14.8 & 1,347.4 & 1.06\% & 0.96\% \\ 
m8n2500s2 & 5,835,360 & 5,895,176 & {\bf 5,887,819} &  13.6 & 1,095.8 & 1.02\%  & 0.89\% \\
m8n2500s3 & 6,046,374 & 6,110,906 & {\bf 6,104,058} &  14.5 & 1.015.1 & 1.06\% & 0.95\% \\
m8n2500s4 & 6,021,410 & 6,078,364 & {\bf 6,075,874} &  14.0 & 1,238.4 & 0.94\% & 0.90\% \\
\hline
m12n1500s0 & 3,621,952 & 3,670,642 & {\bf 3,663,952} & 5.1 & 215.9 & 1.34\% & 1.16\% \\
m12n1500s1 & 3,523,474 & {\bf 3,570,252} & 3,570,484 & 5.5 & 189.0 & 1.32\% & 1.33 \% \\
m12n1500s2 & 3,932,474 & 3,988,062 & {\bf 3,983,324} & 4.6 & 160.4 & 1.41\% & 1.29\% \\
m12n1500s3 & 3,789,274 & 3,833,318 & {\bf 3,831,427} & 4.7 & 122.5 & 1.15\% & 1.10\% \\
m12n1500s4 & 3,694,646 & 3,745,298 & {\bf 3,738,872} & 4.6 & 233.5 & 1.37\% & 1.19\% \\
\hline
m12n2000s0 & 5,239,126 & 5,301,310 & {\bf 5,294,030} & 10.8 & 530.9 & 1.18\% & 1.04\% \\
m12n2000s1 & 4,844,414 & 4,907,954 & {\bf 4,899,575} &  9.4 & 3,275.3 & 1.31\% & 1.13\% \\
m12n2000s2 & 4,611,692 & 4,667,510 & {\bf 4,665,170} &  9.3 & 883.5 & 1.21\% & 1.16\% \\
m12n2000s3 & 4,822,028 & 4,881,702 & {\bf 4,871,930} &  9.1 & 438.5 & 1.23\% & 1.03\% \\
m12n2000s4 & 4,961,406 & 5,025,946 & {\bf 5,019,364} &  8.8 & 684.5 & 1.30\% & 1.16\% \\
\hline
m12n2500s0 & 5,860,766 & 5,948,488 & {\bf 5,937,754} &  14.5 & 1,101.0 & 1.49\% & 1.13\% \\
m12n2500s1 & 6,000,516 & {\bf 6,070,994} & 6,071,772 &  14.4 & 1,534.9 & 1.17\% & 1.18\% \\
m12n2500s2 & 5,940,276 & 6,012,290 & {\bf 6,005,232} &  15.6 & 1,236.5 & 1.21\% & 1.10\% \\
m12n2500s3 & 6,072,130 & 6,154,820 & {\bf 6,140,502} & 14.5 & 1,060.8& 1.36\% & 1.12\% \\
m12n2500s4 & 5,748,976 & 5,817,298 & {\bf 5,819,317} &  15.5 & 1,725.6 & 1.18\% & 1.12\% \\
\hline
m16n1500s0 & 3,568,522 & 3,618,204 & {\bf 3,616,462} & 5.5 & 195.7 & 1.39\% & 1.34\% \\
m16n1500s1 & 3,591,374 & {\bf 3,637,940} & 3,641,353 & 5.5 & 150.2 & 1.29\% & 1.39\% \\
m16n1500s2 & 3,554,800 & 3,604,392 & {\bf 3,601,830} & 4.8 & 322.5 & 1.39\% & 1.32\% \\
m16n1500s3 & 3,861,652 & 3,914,558 & {\bf 3,908,222} & 4.8 & 217.2 & 1.27\% & 1.20\% \\
m16n1500s4 & 3,603,796 & {\bf 3,657,334} & 3,659,326 & 5.6 & 357.9 & 1.48\% & 1.54\%\\
\hline
m16n2000s0 & 4,789,504 & 4,856,390 & {\bf 4,852,389} & 9.8 & 3,705.6 & 1.39\% & 1.31\% \\
m16n2000s1 & 4,680,998 & 4,745,990 & {\bf 4,744,248} & 9.7 & 2,347.0 & 1.38\% & 1.35\% \\
m16n2000s2 & 4,774,408 & 4,847,436 & {\bf 4,834,352} & 10.1 & 732.3 & 1.53\% & 1.25\% \\
m16n2000s3 & 4,850,652 & 4,926,124 & {\bf 4,920,016} & 9.9 & 520.2 & 1.55\% & 1.43\% \\
m16n2000s4 & 4,700,490 & 4,767,256 & {\bf 4,761,708} &  9.6 & 649.6 & 1.42\% &  1.30\%\\
\hline
m16n2500s0 & 5,960,298 & 6,045,500 & {\bf 6,038,849} & 15.5 & 1,318.1 & 1.42\% & 1.31\% \\
m16n2500s1 & 6,055,252 & 6,148,106 & {\bf 6,138,606} & 15.0 & 1,566.7 & 1.53\% & 1.37\% \\
m16n2500s2 & 6,043,364 & 6,123,422 & {\bf 6,118,930} & 16.2 & 1,240.2 & 1.32\% & 1.25\% \\
m16n2500s3 & 6,067,858 & 6,155,328 & {\bf 6,148,939} & 15.3 & 902.3 & 1.44\% & 1.33\% \\
m16n2500s4 & 5,857,966 & 5,952,214 & {\bf 5,942,387} &  14.5 & 1,414.2 & 1.60\% & 1.44\%
\end{tabular}
\end{table}

\section{Concluding remarks}
\label{sect5}

The MDVSP is very important in the management of public transport systems. Our paper introduces three heuristics based on a classical integer linear programming formulation and on graph theoretic methods of fixing the infeasible subtours gathered from an integer solution.

The effectiveness of our heuristics was proved by extensive experimentations using a large set of benchmark instances. Our heuristics are fast and give good results compared with other existing solutions; the running times are much lower (but we performed the experiments four years later) and the results are better than in \cite{guedes16} for the benchmarks from Huisman's website. Compared with the results in \cite{kulkarni18}, the running times for larger benchmark instances are lower and the results remain within $ 1.6\% $ percent error, while we worked on a mainstream desktop PC instead of a dedicated server.

Future work will be directed towards a truncated branch and bound (and branch and cut) algorithm based on adding constraints like in the second relaxation - corresponding with a column generation procedure for the dual problem. Our heuristics can be used for providing initial upper bounds for such an algorithm.

\end{document}